\documentclass{article}
\usepackage{epsfig}

\begin{document}
\title{Numerical Solution of Differential Equations in Irregular Plane Regions Using Quality Structured Convex Grids}

\author{F. Dom\'inguez-Mota, M. Equihua, S. Mendoza\\ and J.~G.~Tinoco-Ruiz\\\\
\em Facultad de Ciencias F\'isico Matem\'aticas,\\
\em Universidad Michoacana de San Nicol\'as de Hidalgo.\\
\em Edificio ``B", Ciudad Universitaria, Morelia, C.P. 58040\\
\em dmota@umich.mx
}
\renewcommand{\today}{October 6, 2011}
\maketitle

\begin{abstract}
The variational grid generation method is a powerful tool for generating structured convex grids on irregular simply connected domains whose boundary is a polygonal Jordan curve. Several examples that show the accuracy of a difference approximation to the solution of a Poisson equation  using these kind of structured grids have been recently reported. In this paper, we compare the accuracy of the numerical solution calculated by applying those structured grids with finite differences against the the solution obtained with Delaunay-like triangulations on irregular regions.
\end{abstract}


\section{Introduction}
\begin{figure}[htb]
\begin{center}
\epsfig{file=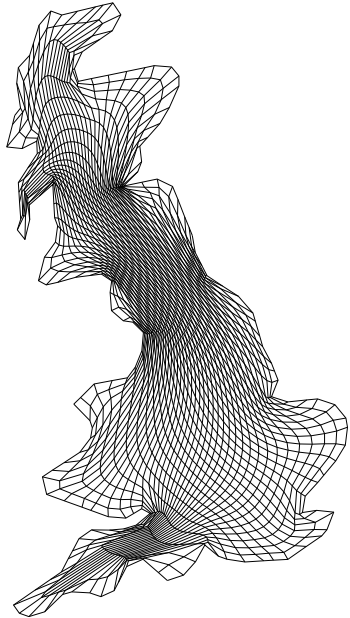, width=1.1\textwidth}
\end{center}
\caption{Structured grid generated by UNAMALLA}
\label{suave}
\end{figure}

\begin{figure}[htb]
\begin{center}
\epsfig{file=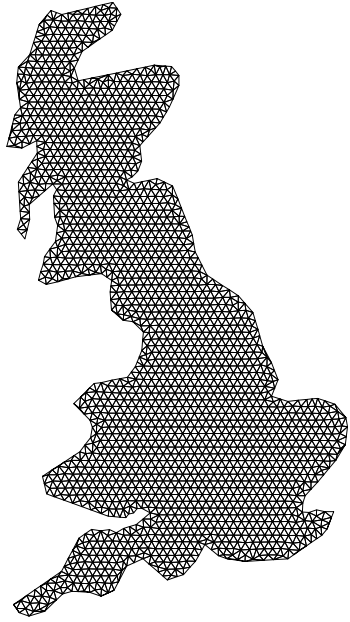, width=1.1\textwidth}
\end{center}
\caption{Triangulation generated by DistMesh}
\label{rrr42}
\end{figure}

For the numerical solution of the Poisson equation on irregular domains, the use of finite differences and finite elements with triangulations obtained by subdividing each grid cell of a structured grid along a diagonal has been addressed in Ref.~\cite{gto}, proving that the finite difference approach on such grids is accurate enough. However, since structured grids often have some elongated cells, a natural question that arises in this context is whether this numerical solution is more accurate than the solution obtained by using finite elements on a standard Delaunay-like triangulation.\\
In this paper, we compare the accuracy of the numerical solution for this problem using finite differences in structured grids, and finite elements on Delaunay-like triangulations. We are specifically interested in irregular boundaries, since their geometry is closer to a realistic domain, for instance, a lake.\\
In order to generate the structured convex grids, a variational method was used. The latter consists of minimizing an appropriate functional \cite{castillo2,ivan4,knupp0,steinbergroache}. Area and harmonic
functionals can be used for gridding a wide variety of  simple connected domains in the plane \cite{harmonic,morelia,motabarrera,knuppsteinberg,tincoco,tinoco3}, whose boundaries are
closed polygonal Jordan curves with positive orientation. \\
If  $m$ and $n$ represent the ``vertical" and ``horizontal" numbers of points of the ``sides", then the boundary is the positively oriented polygonal curve $\gamma$ of vertices
$V~=~\{v_1,\cdots,v_{2(m+n-2)}\}$, and it defines the typical domain $\Omega$.\\
A doubly indexed set \begin{displaymath}G=\{P_{i,j}|1\le i\le~m, 1\le j\le n\}\end{displaymath} of points of the plane
with the fixed boundary positions given by $V$ is a logically rectangular structured grid with quadrilateral elements for $\Omega$, of order $m\times n$. \\
A grid $G$ is convex if and only if each one of the $(m-1)(n-1)$ quadrilaterals (or cells) $c_{i,j}$ of vertices
$\{P_{i,j},P_{i+1,j},P_{i,j+1},P_{i+1,j+1}\}$, $1\le i<m$, $1\le j<n$, is convex and non-degenerate (See fig. \ref{suave}).\\
The basis for the direct optimization method, as developed by Ivanenko {\em et al.} \cite{ivan4}, is the minimization a suitable function of the form
\begin{equation}\label{discreto}
F(G)=\sum_{i=1}^{m-1}\sum_{j=1}^{n-1} f(c_{i,j}),
\end{equation}
where $c_{i,j}$ is the $(i,j)^{th}$ grid cell  and $f$ is a function of its vertices; the problem is to find the coordinates of the interior points of the grid $G$. The functional used to generate the structured grids of the numerical tests, as implemented in UNAMALLA \cite{unamalla}, was the adaptive linear convex combination of the area functional $S_{\omega}(G)$ described in Ref.~\cite{etna}, and the length functional $L(G)$ with weight $\sigma=0.5$ (See Ref.~\cite{harmonic}).\\
The parameter $\omega$, which is a scale factor, can be updated in such a way that in a finite number of updates
the combined functional attain its minima within the set of convex grids for $\Omega$ if the latter is nonempty.
Further properties of this functional, as well as the algorithm for updating its parameter has been reported in Ref.~\cite{harmonic} and Ref.~\cite{etna}.\\
The triangulations we considered for this paper are those generated by DistMesh \cite{persson}, which is based on the physical analogy between a simplex mesh and a truss structure, where the meshpoints are nodes of the truss. It generates an initial Delaunay triangulation, then assumes an appropriate force-displacement function for the bars in the truss at each iteration, and finally solves for equilibrium (See fig. \ref{rrr42}). In order to have grids with a similar number of elements, we used two initializations: a) the initial edge length for DistMesh was set in proportion to half the average diagonal length of the structured grids, b) a variation of DistMesh was designed for which the initial points inside the region are the inner nodes of the corresponding structured grid.\\
It must be noted that for very irregular boundaries,  DistMesh might produce a few triangular elements along the boundary which do not satisfy the Delaunay condition.\\
\section{Finite difference approximation}\label{fda}
Standard difference schemes can be generalized by considering a finite set of nodes $p_1,p_2,...,p_k$, for which it is required to find
coefficients $\Gamma_1,\Gamma_2,...,\Gamma_k$ such that \cite{celia}
\begin{equation}\label{defdf}
\frac{\partial^qu}{\partial x^{l} \partial y^{q-l}}\vert_{x=x_{*}}\approx \sum_{i}\Gamma_i u(p_i).
\end{equation}
As it is well known, the $\Gamma$ values can be calculated with ease in regular regions. However, despite the basic idea is quite simple, the application to Taylor's Theorem leads to more complicated schemes on irregular regions; up to our knowledge, there are few efficient schemes for such kind of regions.\\
The calculation of these coefficients has been studied by Tinoco {\it et al} \cite{azucena}, and \
Shashkov \cite{shashkov}.
We make use of the second order scheme developed for the boundary value problem
\begin{eqnarray}\label{eqfuerte}
-\nabla(K(x,y)\nabla u(x,y))&=&f(x,y)\\
 K(x,y)&=&\left( \begin{array}{rrr} K_{1 1} & & K_{1 2} \\ K_{1 2} & & K_{2 2} \end{array} \right)\nonumber\\
u(x,y)\vert_{\partial \Omega}&=&g(x,y),\nonumber
\end{eqnarray} with non-diagonal matrices $K(x,y)$ studied in Ref.~\cite{shashkov}, which is based on the method of support-operators and has the advantage of providing explicit formulas for the $\Gamma$ coefficients\footnote{This scheme is second order according to the grid norm defined by equation (\ref{normis}).} .\\
In the numerical experiments, we selected the $3\times 3$ subgrid defined for the nodes $x(i-1:i+1,j-1:j+1)$, $y(i-1:i+1,j-1:j+1)$ to approximate $-\nabla(K(x,y)\nabla u(x,y))=f(x,y)$ at the inner grid node $(x(i,j),y(i,j))$. As in the rectangular case, an algebraic system of equations is obtained from discretization, which becomes sparse as $m$ and $n$ increase.\\
\section{Finite elements approximation}\label{fea}
Let $Ne$ be the number of triangular elements in a grid.  Galerkin's approach to the solution of equation (\ref{eqfuerte}) is given
by the combination
\begin{equation}\label{defu}
 u(x,y)\approx \sum_{i=1}^{Ne}U_i \phi_i
\end{equation}
of trial-test pyramid functions whose faces are defined on a triangle by
\begin{equation}\label{le}
\phi(x,y)=A+Bx+Cy.
\end{equation}
This selection yields the weak formulation
\begin{equation}\label{deffem}
\sum_{i=1}^{Ne} U_i\int_{\Omega}<\nabla\phi_j,K(x,y)\nabla \phi_i>d\mathcal{A}=\int_{\Omega}\phi_j f d\mathcal{A},~~~j=1,\cdots,Ne,
\end{equation}
where $<\cdot,\cdot>$ is the canonical inner product (See Ref.~\cite{strang}).\\
There are several other possible choices of elements. However, as it is well known, the use of pyramids on triangulations allows a very simple second order approximation to the solution of equation (\ref{eqfuerte}).\\
\section{Numerical tests}
For the numerical tests, we have selected 9 polygonal regions, most of them approximations to real geographical locations: they will be denoted as Dome (dom), Great Britain (eng), Havana bay (hab), M19 (m19), M\'exico (mex), Plow (plo), Swan (swa), Ucha (uch) and Michoac\'an (mic).
They are shown in figure 3. \\
\begin{figure}[htb]
\hbox to \hsize {\hfill
\begin{minipage}{3.3cm}
\begin{center}
\epsfig{file=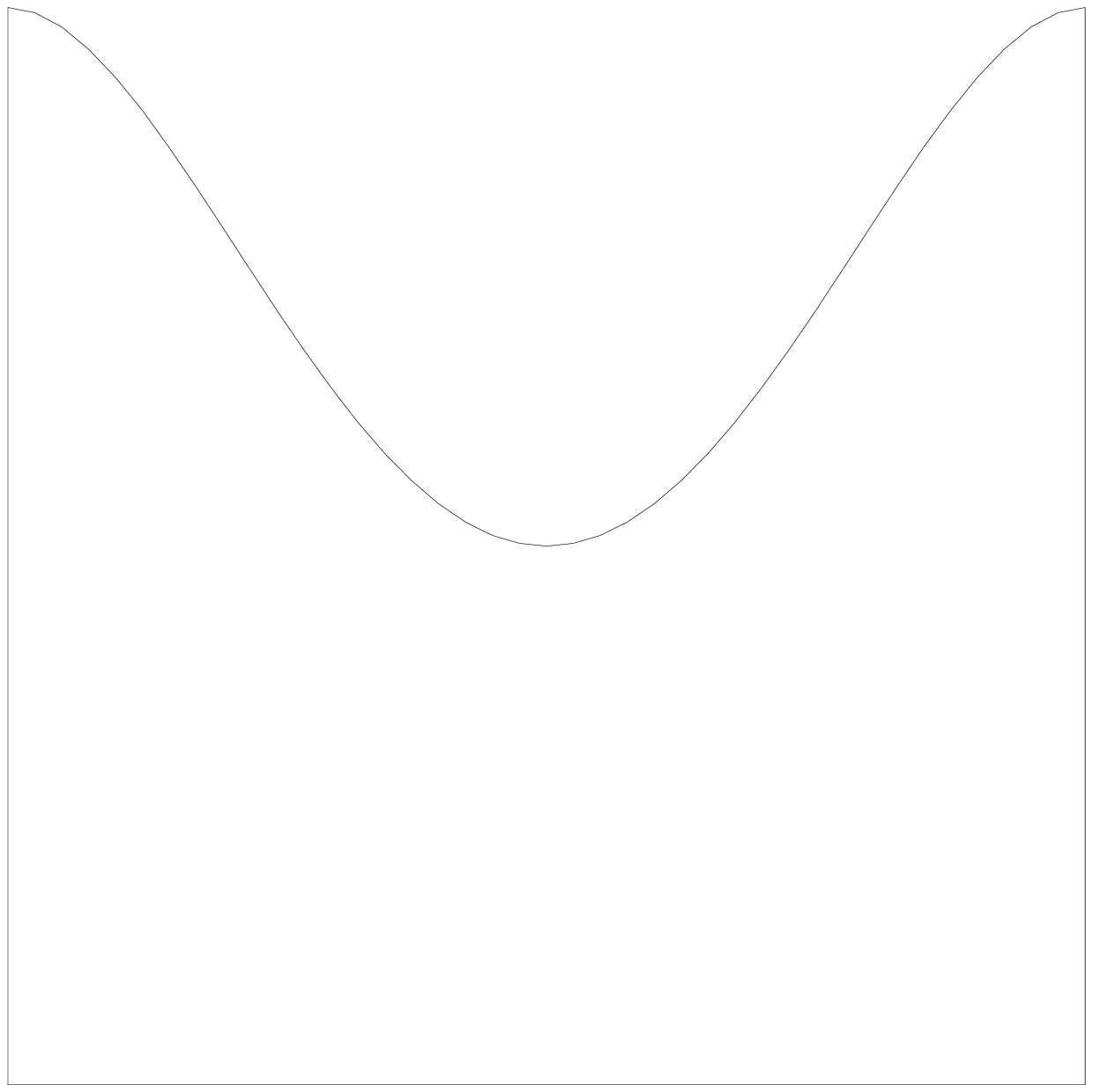, width=3cm, height=3cm}\\
Dome\\
\epsfig{file=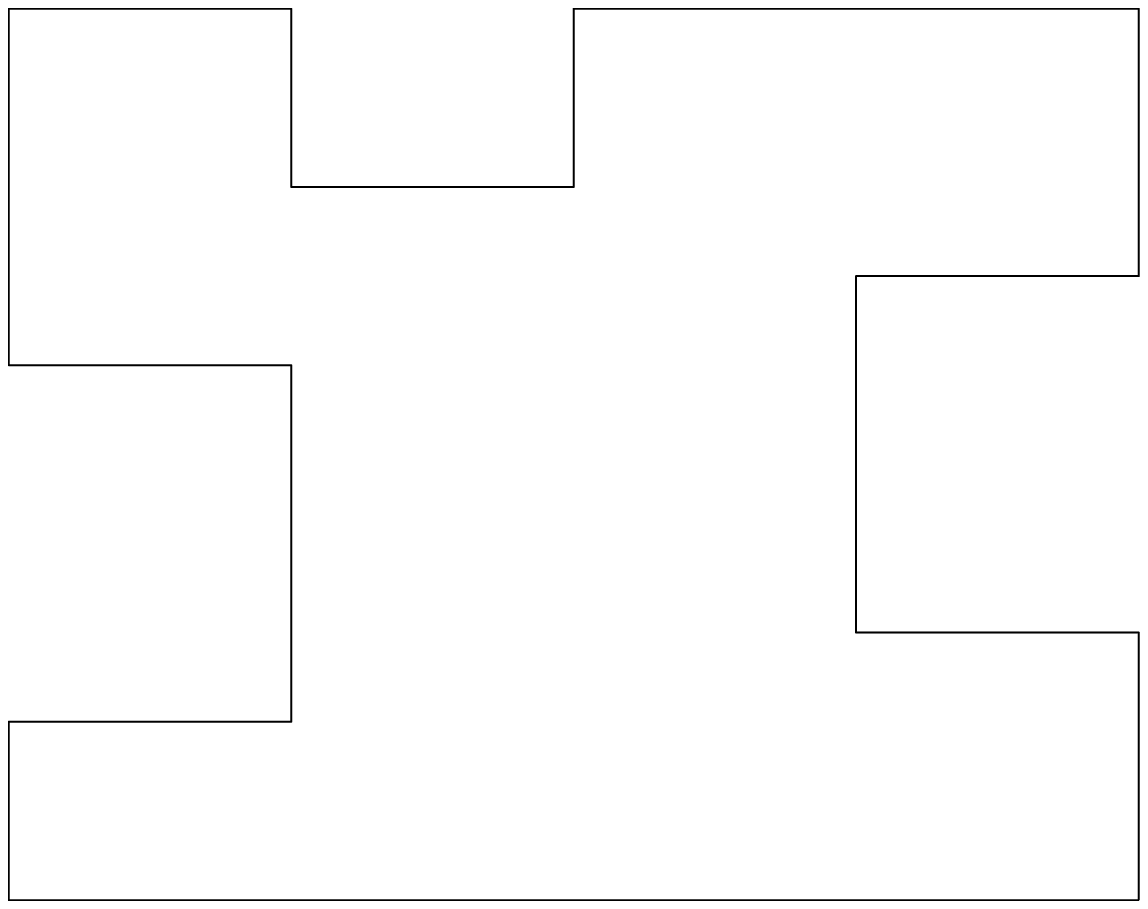, width=3cm, height=3cm}\\
M19\\
\epsfig{file=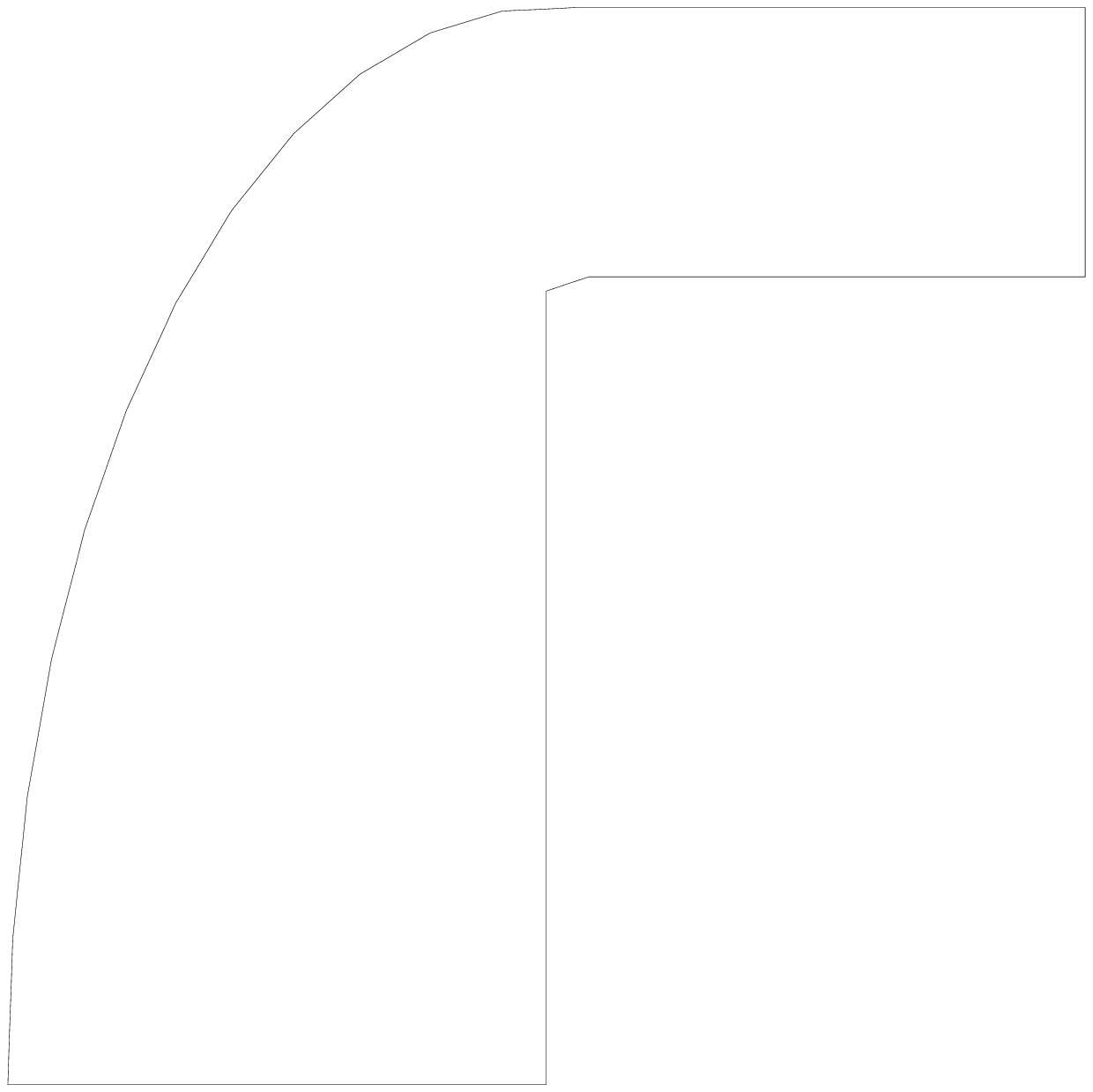, width=3cm, height=3cm}\\
Plow\\
\end{center}
\end{minipage}
\hfill
\begin{minipage}{3.3cm}
\begin{center}
\epsfig{file=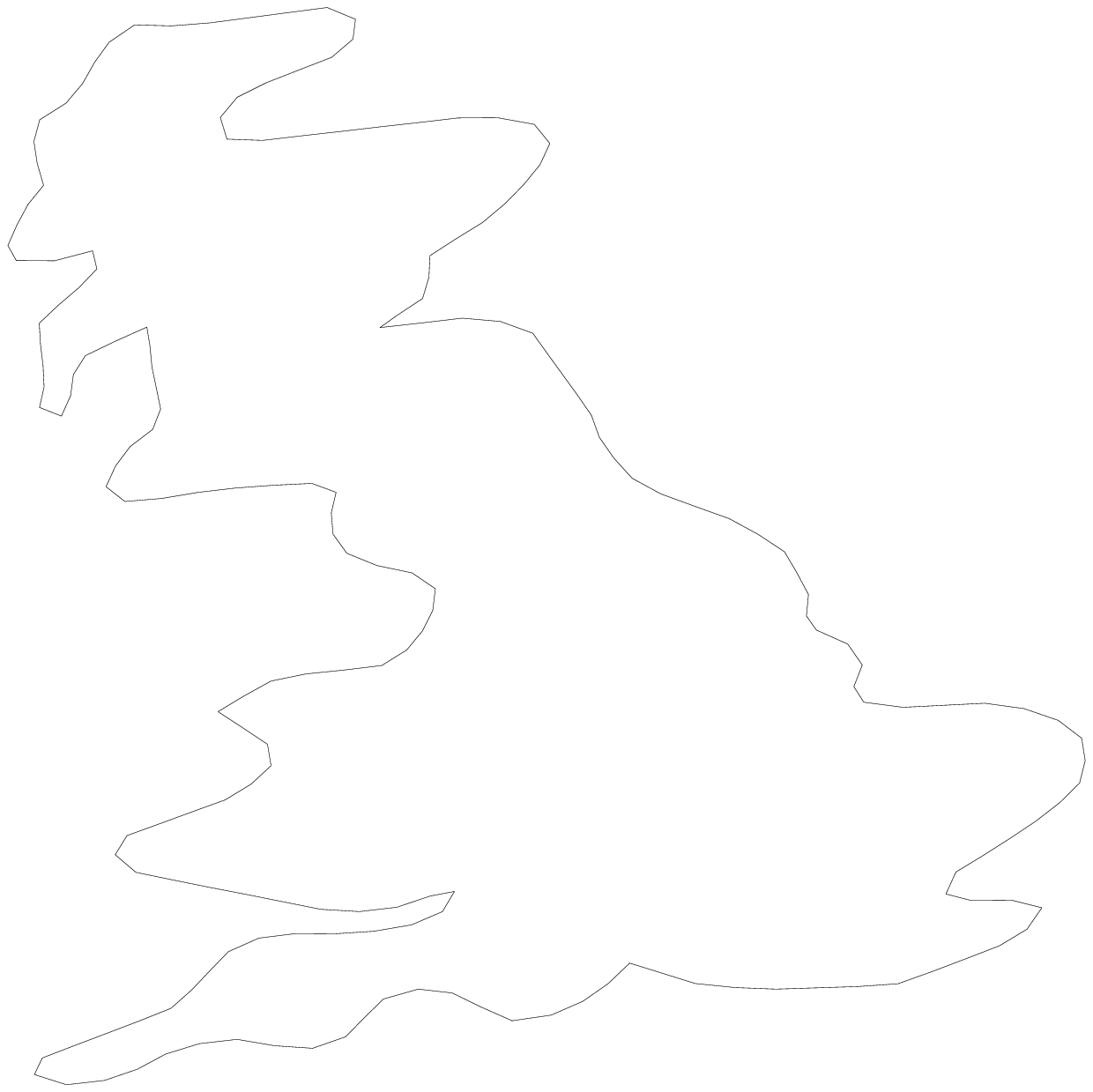, width=3cm, height=3cm}\\
Great Britain\\
\epsfig{file=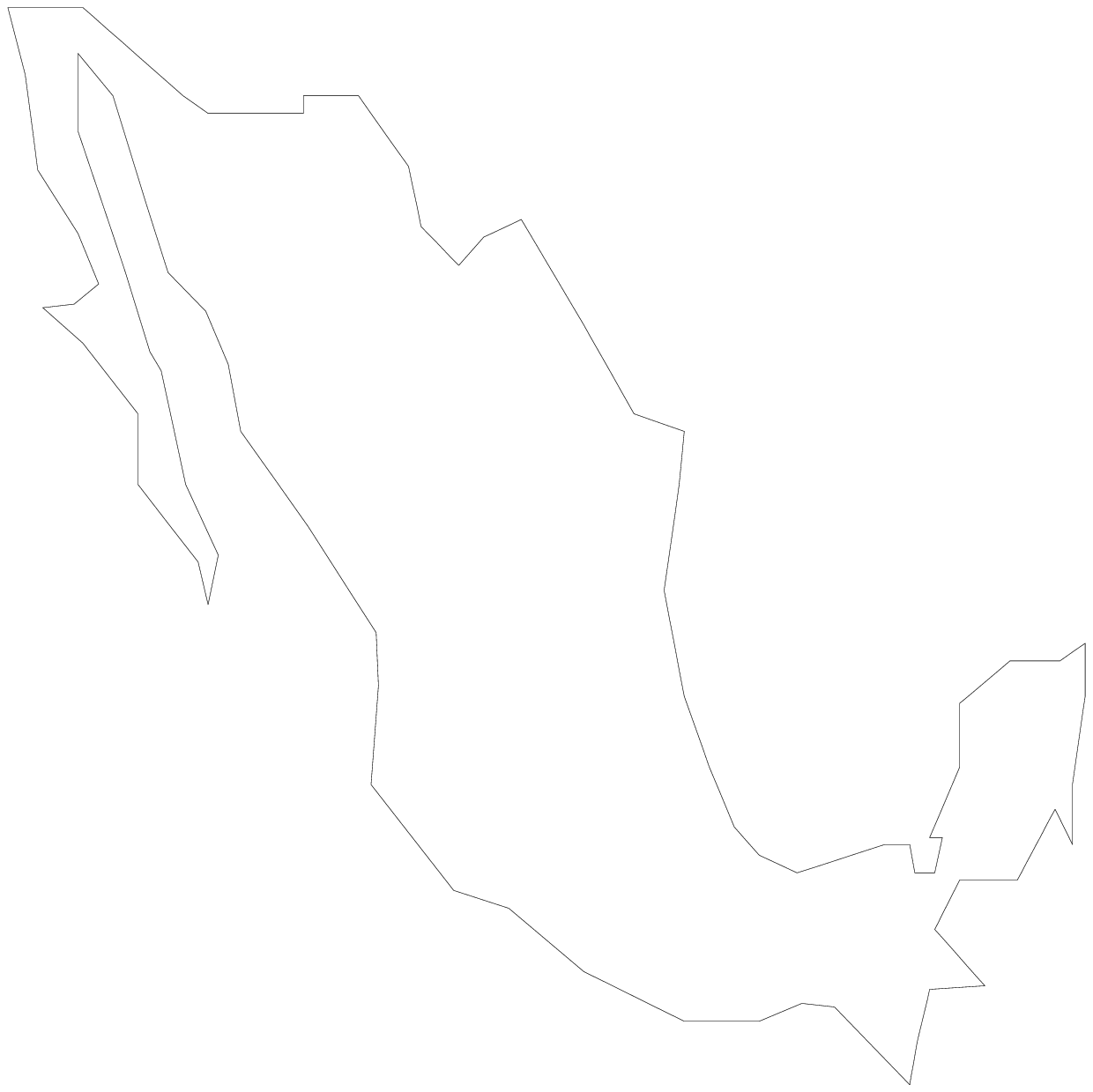, width=3cm, height=3cm}\\
M\'exico\\
\epsfig{file=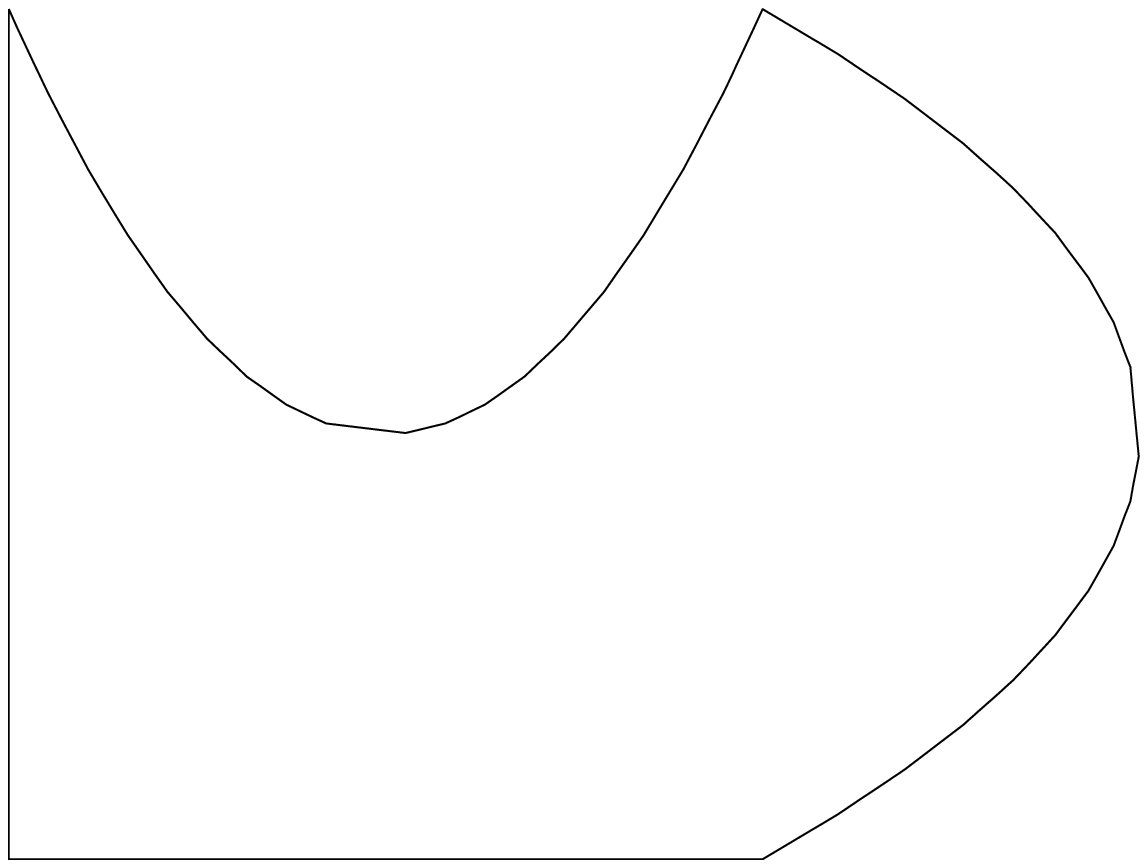, width=3cm, height=3cm}\\
Swan\\
\end{center}
\end{minipage}
\hfill
\begin{minipage}{3.3cm}
\begin{center}
\epsfig{file=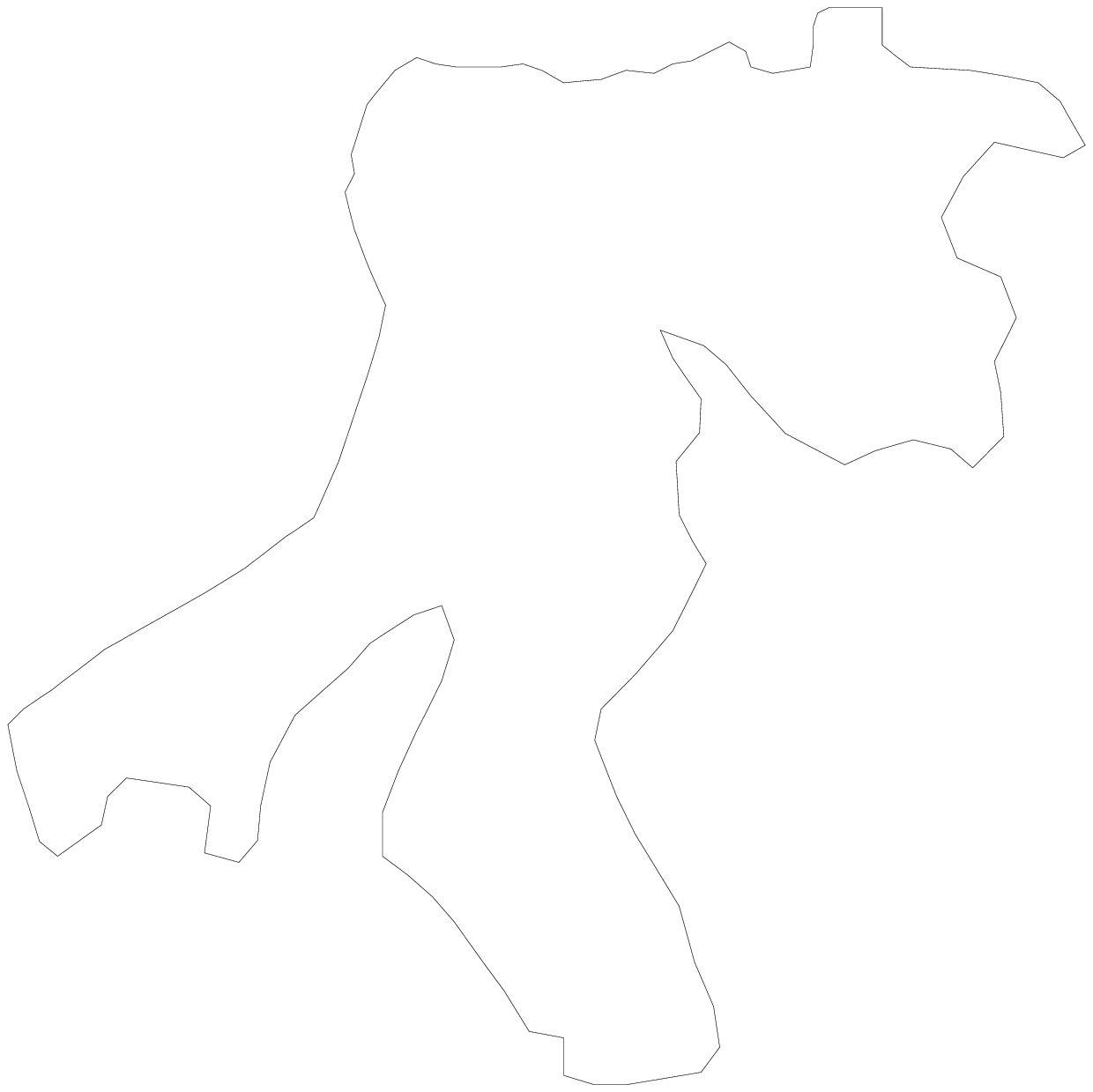, width=3cm, height=3cm}\\
Havana bay\\
\epsfig{file=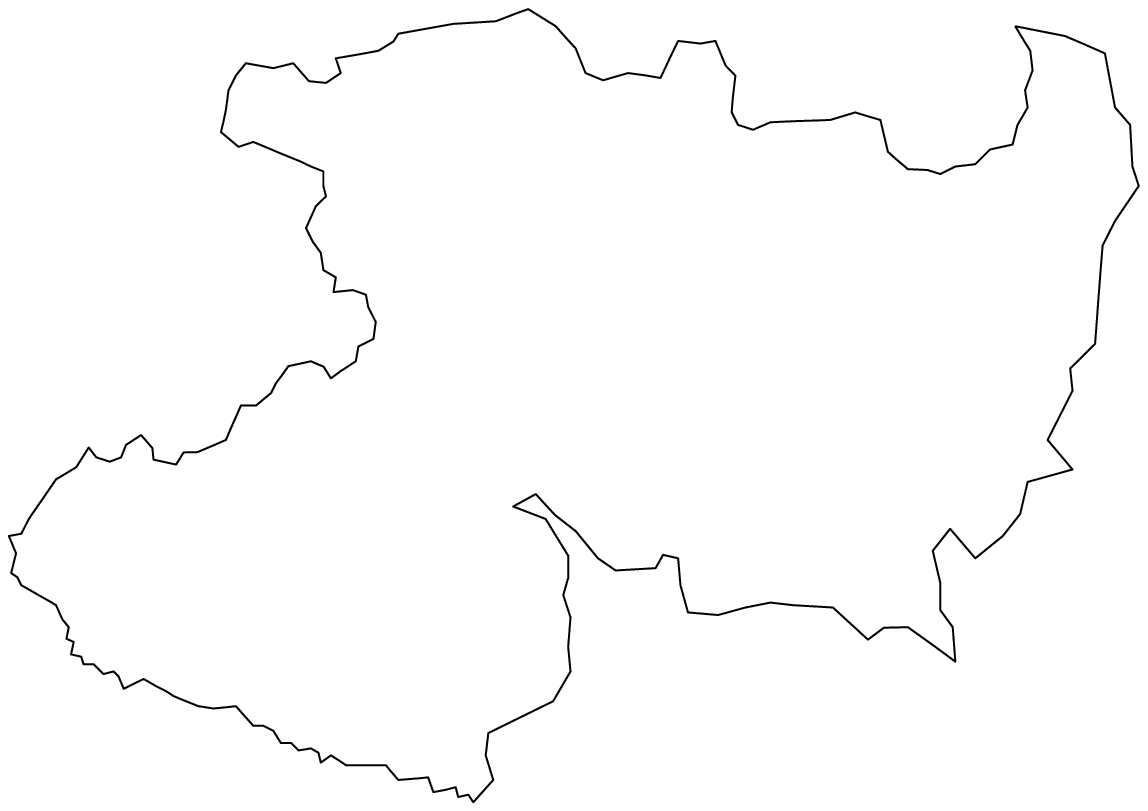, width=3cm, height=3cm}\\
Michoac\'an\\
\epsfig{file=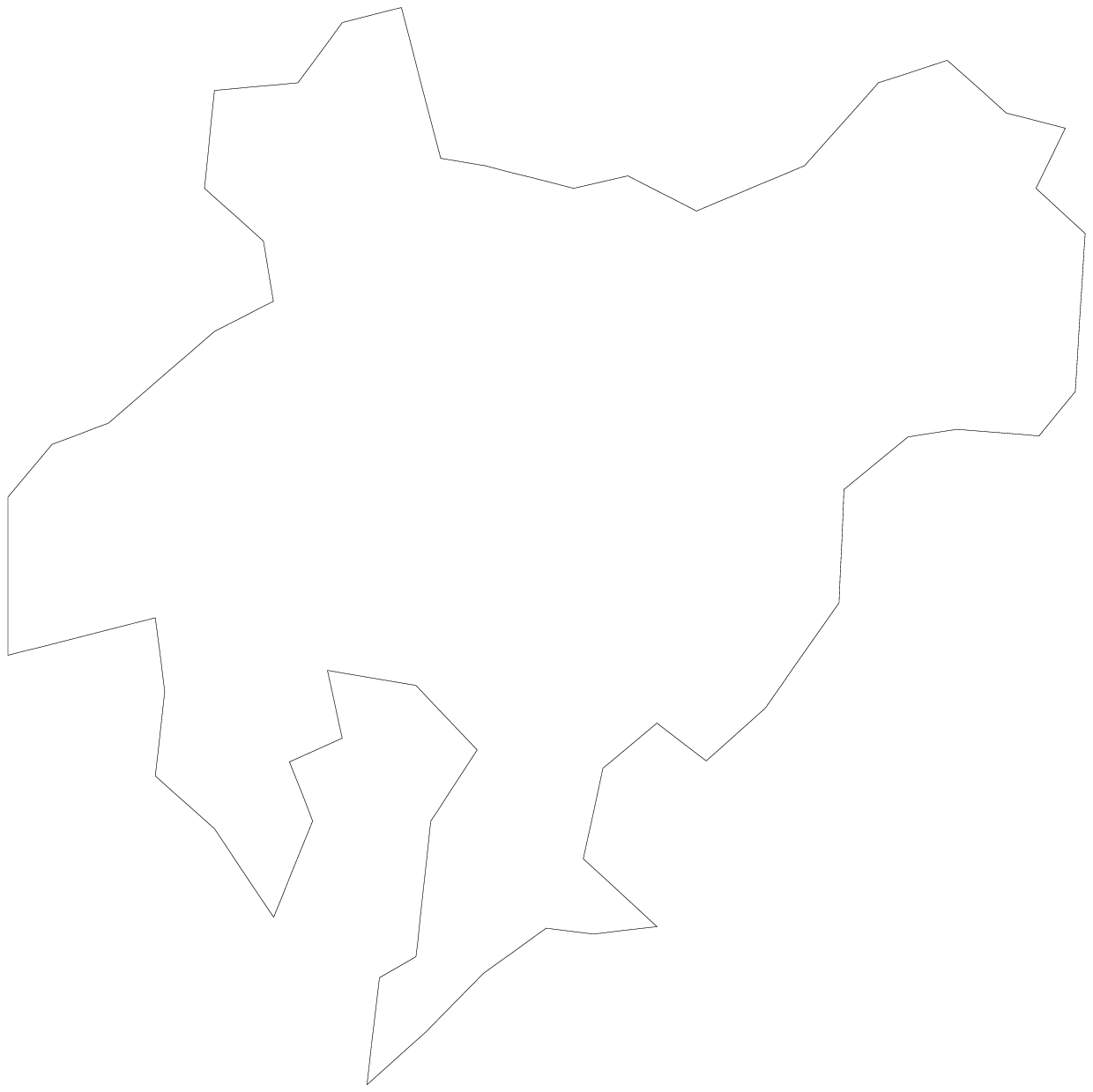, width=3cm, height=3cm}\\
Ucha\\
\end{center}
\end{minipage}
\hfill}
\caption{Test regions.}
\label{B}
\end{figure}
\noindent Scaling these boundaries in order to lie in $[0,1]\times[0,1]$, convex grids with 21, 41 and 81 points per side were generated by minimizing the $1/2(S_\omega(G)+L(G))$  functional in UNAMALLA with default parameters. The resulting structured grids were used with Shashkov's finite difference schemes \cite{shashkov}. As mentioned before, they were also used as initial data for some triangulations. \\
DistMesh was used to triangulate the same test polygonal boundaries with default parameters, setting the bar length as half the average diagonal length calculated in the corresponding structured grids (denoted as $DistMesh_a$). A variation of DistMesh was also designed, for which the initial points inside the region are the inner nodes of the corresponding structured grid: these grids are denoted as $DistMesh_b$ grids. \\
The number of elements and inner nodes in each grid can be seen in in table \ref{tab:ne}. The column $N$ gives twice the number of grid cells\footnote{This number is considered because the number of triangles in a triangulation obtained by subdividing each grid cell along a diagonal of a structured grid is twice the number of cells.},  and the columns $N_a$ and $N_b$ the number of triangular elements in the $DistMesh_a$ and $DistMesh_b$ triangulations, respectively. The columns labeled $Nu$, $Nu_a$ and $Nu_b$ represent the corresponding number of inner nodes, which is equal to the number of unknowns in the approximation.\\
\begin{table}
\centering
\caption{Number of grid elements and inner nodes}
\label{tab:ne}       
\begin{tabular}{lrrrrrrr}
\hline\noalign{\smallskip}
\noalign{\smallskip}
Region	&	Size	&	$N_1$	&	$N_2$&$N_3$	&	$Nu$ & $Nu_a$	& $Nu_b$\\
\noalign{\smallskip}\hline\noalign{\smallskip}
dom	&	21	&	800	&	1540	&	795	&	361	&	693	&	356	\\
	&	41	&	3200	&	6064	&	3189	&	1521	&	2884	&	1509	\\
	&	81	&	12800	&	24150	&	12783	&	6241	&	11785	&	6218	\\
eng	&	21	&	800	&	841	&	771	&	361	&	355	&	337	\\
	&	41	&	3200	&	3463	&	3111	&	1521	&	1581	&	1440	\\
	&	81	&	12800	&	13918	&	12614	&	6241	&	6622	&	6074	\\
hab	&	21	&	800	&	1261	&	767	&	361	&	559	&	333	\\
	&	41	&	3200	&	5048	&	3115	&	1521	&	2370	&	1446	\\
	&	81	&	12800	&	20313	&	12605	&	6241	&	9845	&	6070	\\
m19	&	21	&	800	&	1420	&	785	&	361	&	626	&	350	\\
	&	41	&	3200	&	5717	&	3157	&	1521	&	2696	&	1483	\\
	&	81	&	12800	&	22772	&	12705	&	6241	&	11067	&	6156	\\
mex	&	21	&	800	&	477	&	686	&	361	&	183	&	273	\\
	&	41	&	3200	&	1781	&	2911	&	1521	&	771	&	1288	\\
	&	81	&	12800	&	7036	&	11918	&	6241	&	3301	&	5519	\\
mic	&	21	&	800	&	1432	&	796	&	361	&	645	&	356	\\
	&	41	&	3200	&	5760	&	3183	&	1521	&	2747	&	1506	\\
	&	81	&	12800	&	22812	&	12694	&	6241	&	11096	&	6136	\\
plo	&	21	&	800	&	863	&	748	&	361	&	363	&	310	\\
	&	41	&	3200	&	3335	&	3076	&	1521	&	1535	&	1399	\\
	&	81	&	12800	&	13115	&	12537	&	6241	&	6303	&	5983	\\
swa	&	21	&	800	&	1410	&	796	&	361	&	634	&	357	\\
	&	41	&	3200	&	5495	&	3187	&	1521	&	2612	&	1508	\\
	&	81	&	12800	&	21737	&	12766	&	6241	&	10580	&	6205	\\
uch	&	21	&	800	&	1088	&	794	&	361	&	483	&	355	\\
	&	41	&	3200	&	4252	&	3166	&	1521	&	1976	&	1489	\\
	&	81	&	12800	&	17160	&	12729	&	6241	&	8288	&	6172	\\
\noalign{\smallskip}\hline
\end{tabular}
\end{table}

The resulting algebraic systems for finite differences and elements are sparse due to the discretization. However, the systems for finite differences are block-tridiagonal, so they can be solved with a number of algorithms in a very efficient way; this is a clear advantage of the double index in a structured grid. For the tests, Gauss-Seidel Method was used.\\
On the other hand, the matrices in the finite element systems have no specific structure and require a more complicated data structure in order to solve them efficiently. For the tests, they were solved using a sparse Gaussian elimination routine.\\
The following values for $u$ and $K$ were selected (See Ref.~\cite{shashkov}):
\begin{enumerate}
\item First problem. $$K(x,y)=\left( \begin{array}{rrr} 1 & & 0 \\ 0 & & 1 \end{array} \right),\quad u=2\exp(2x+y).$$
\item Second problem. $$K(x,y)=P^TDP,$$
with $$P=\left( \begin{array}{rrr} cos(\frac{\pi}{8}) & & sin(\frac{\pi}{8}) \\ -sin(\frac{\pi}{8}) & & cos(\frac{\pi}{8}) \end{array} \right)$$ and $$D=\left( \begin{array}{rrr} 1+2x^2+y^2 & & 0 \\ 0 & & 1+x^2+2y^2 \end{array} \right)$$
$$u=\sin(\pi x)\sin(\pi y).$$
\item Third problem. $$K(x,y)=P^TDP,$$
with $$P=\left( \begin{array}{rrr} cos(\frac{\pi}{4}) & & sin(\frac{\pi}{4}) \\ -sin(\frac{\pi}{4}) & & cos(\frac{\pi}{4}) \end{array} \right)$$ and $$D=\left( \begin{array}{rrr} 1+2x^2+y^2+y^5 & & 0 \\ 0 & & 1+x^2+2y^2+x^3 \end{array} \right)$$
$$u=\sin(\pi x)\sin(\pi y).$$
\end{enumerate}
Function $f$ was chosen in such a way that $u$ was the exact solution in every case.\\
The $\|\cdot\|_{2}$ error norm for the tests is summarized in tables \ref{tab:1}, \ref{tab:2} and \ref{tab:3}; it was calculated as a grid function
\begin{equation}\label{normis}
\|u-U\|_2=\sqrt{\sum_{i}(u_{i}-U_{i})^2\mathcal{A}_{i}} \quad,
\end{equation}
where $u$ and $U$ are the approximated and the exact solution calculated at the $i^{th}$-element, and $\mathcal{A}_{i}$ is the area of the element. The approximation corresponding to label $Structured$ are those of the second order finite differences mentioned in section \ref{fda}; the approximations for the columns labeled $DistMesh_a$, and $DistMesh_b$ were calculated with finite elements using the pyramid trial-test approximation described in section \ref{fea}. The empirical orders $O$, $O_a$ and $O_b$ between two consecutive grid orders were calculated according to the formula
\begin{equation}
\log\left(\frac{E_i}{E_j}\right)/\log\left(\frac{n_j}{n_i}\right)
\end{equation}
where $E_i$ is the quadratic error associated to the numerical solution calculated with a grid with $n_i$ points per side.
\\
\begin{table}
\centering
\caption{Quadratic error for problem 1}
\label{tab:1}       
\begin{tabular}{lrrrrrrr}
\hline\noalign{\smallskip}
\noalign{\smallskip}
Region	&	Size	&	$Structured$	& $O$ &	$DistMesh_a$	&	$O_1$ & $DistMesh_b$ & $O_2$	\\
\noalign{\smallskip}\hline\noalign{\smallskip}
dom	&	21	&	4.59E-03	&		&	5.77E-03	&		&	6.25E-03	&		\\
	&	41	&	1.22E-03	&	1.98	&	8.79E-04	&	2.81	&	1.47E-03	&	2.16	\\
	&	81	&	1.46E-04	&	3.12	&	1.51E-04	&	2.58	&	2.70E-04	&	2.49	\\
eng	&	21	&	2.58E-03	&		&	1.72E-03	&		&	2.44E-03	&		\\
	&	41	&	4.27E-04	&	2.69	&	6.95E-05	&	4.79	&	4.75E-04	&	2.45	\\
	&	81	&	1.17E-04	&	1.90	&	1.03E-04	&	-0.57	&	2.86E-04	&	0.75	\\
hab	&	21	&	7.53E-03	&		&	1.26E-03	&		&	1.98E-03	&		\\
	&	41	&	1.53E-03	&	2.38	&	1.69E-04	&	3.00	&	5.85E-04	&	1.82	\\
	&	81	&	3.97E-04	&	1.99	&	6.11E-05	&	1.50	&	1.77E-04	&	1.76	\\
m19	&	21	&	5.90E-03	&		&	2.02E-03	&		&	1.94E-02	&		\\
	&	41	&	1.44E-03	&	2.11	&	6.55E-04	&	1.68	&	7.18E-04	&	4.93	\\
	&	81	&	3.21E-04	&	2.20	&	8.24E-04	&	-0.34	&	1.39E-04	&	2.41	\\
mex	&	21	&	8.67E-03	&		&	2.06E-03	&		&	8.03E-04	&		\\
	&	41	&	1.87E-03	&	2.29	&	2.34E-04	&	3.25	&	3.72E-04	&	1.15	\\
	&	81	&	3.85E-04	&	2.32	&	4.93E-05	&	2.29	&	3.82E-05	&	3.34	\\
mic	&	21	&	5.53E-03	&		&	1.45E-04	&		&	8.15E-04	&		\\
	&	41	&	1.77E-03	&	1.70	&	7.00E-05	&	1.09	&	9.21E-04	&	-0.18	\\
	&	81	&	3.73E-04	&	2.29	&	3.45E-04	&	-2.34	&	4.61E-05	&	4.40	\\
plo	&	21	&	1.04E-03	&		&	1.52E-03	&		&	4.49E-04	&		\\
	&	41	&	5.06E-04	&	1.08	&	6.05E-04	&	1.38	&	1.00E-04	&	2.24	\\
	&	81	&	7.97E-05	&	2.71	&	1.77E-04	&	1.81	&	2.80E-05	&	1.87	\\
swa	&	21	&	2.57E-03	&		&	1.08E-03	&		&	1.77E-03	&		\\
	&	41	&	6.01E-04	&	2.17	&	3.73E-04	&	1.58	&	4.54E-04	&	2.03	\\
	&	81	&	1.50E-04	&	2.04	&	5.45E-05	&	2.82	&	9.28E-05	&	2.33	\\
uch	&	21	&	1.14E-02	&		&	3.14E-03	&		&	3.05E-03	&		\\
	&	41	&	1.91E-03	&	2.67	&	1.08E-04	&	5.04	&	7.05E-04	&	2.19	\\
	&	81	&	3.78E-04	&	2.38	&	7.68E-04	&	-2.88	&	1.84E-04	&	1.97	\\
\noalign{\smallskip}\hline
\end{tabular}
\end{table}

\begin{table}
\centering
\caption{Quadratic error for problem 2}
\label{tab:2}       
\begin{tabular}{lrrrrrrr}
\hline\noalign{\smallskip}
\noalign{\smallskip}
Region	&	Size	&	$Structured$	& $O$ &	$DistMesh_a$	&	$O_1$ & $DistMesh_b$ & $O_2$	\\
\noalign{\smallskip}\hline\noalign{\smallskip}
dom	&	21	&	9.68E-04	&		&	4.21E-04	&		&	5.40E-04	&		\\
	&	41	&	2.26E-04	&	2.18	&	7.54E-05	&	2.57	&	1.44E-04	&	1.98	\\
	&	81	&	4.75E-05	&	2.29	&	1.04E-05	&	2.91	&	3.16E-05	&	2.23	\\
eng	&	21	&	9.72E-04	&		&	6.31E-04	&		&	7.68E-04	&		\\
	&	41	&	1.83E-04	&	2.50	&	2.49E-05	&	4.83	&	2.63E-04	&	1.60	\\
	&	81	&	5.00E-05	&	1.90	&	4.40E-05	&	-0.84	&	1.30E-04	&	1.04	\\
hab	&	21	&	1.08E-03	&		&	1.95E-04	&		&	3.54E-04	&		\\
	&	41	&	3.90E-04	&	1.52	&	2.79E-05	&	2.91	&	6.58E-05	&	2.51	\\
	&	81	&	9.72E-05	&	2.04	&	1.37E-05	&	1.04	&	1.99E-05	&	1.76	\\
m19	&	21	&	1.30E-03	&		&	2.93E-04	&		&	2.34E-03	&		\\
	&	41	&	2.75E-04	&	2.32	&	2.52E-05	&	3.66	&	1.26E-04	&	4.37	\\
	&	81	&	8.46E-05	&	1.73	&	1.33E-04	&	-2.45	&	2.91E-05	&	2.15	\\
mex	&	21	&	1.96E-03	&		&	5.39E-04	&		&	3.24E-04	&		\\
	&	41	&	3.00E-04	&	2.80	&	8.17E-05	&	2.82	&	9.53E-05	&	1.83	\\
	&	81	&	7.54E-05	&	2.03	&	1.46E-05	&	2.53	&	4.65E-06	&	4.44	\\
mic	&	21	&	8.01E-04	&		&	2.78E-05	&		&	1.08E-04	&		\\
	&	41	&	2.03E-04	&	2.05	&	5.99E-05	&	-1.15	&	7.51E-05	&	0.55	\\
	&	81	&	6.72E-05	&	1.62	&	3.90E-05	&	0.63	&	4.65E-06	&	4.09	\\
plo	&	21	&	2.36E-04	&		&	2.87E-04	&		&	1.13E-04	&		\\
	&	41	&	9.72E-05	&	1.33	&	1.18E-04	&	1.33	&	2.26E-05	&	2.41	\\
	&	81	&	2.28E-05	&	2.13	&	2.31E-05	&	2.40	&	4.98E-06	&	2.22	\\
swa	&	21	&	4.63E-04	&		&	2.35E-04	&		&	3.66E-04	&		\\
	&	41	&	1.14E-04	&	2.10	&	7.44E-05	&	1.72	&	8.81E-05	&	2.13	\\
	&	81	&	3.37E-05	&	1.79	&	1.15E-05	&	2.75	&	2.04E-05	&	2.15	\\
uch	&	21	&	1.33E-03	&		&	9.61E-04	&		&	4.75E-04	&		\\
	&	41	&	2.59E-04	&	2.44	&	2.29E-05	&	5.59	&	1.30E-04	&	1.94	\\
	&	81	&	6.14E-05	&	2.12	&	1.33E-04	&	-2.58	&	3.00E-05	&	2.15	\\
\noalign{\smallskip}\hline
\end{tabular}
\end{table}

\begin{table}
\centering
\caption{Quadratic error for problem 3}
\label{tab:3}       
\begin{tabular}{lrrrrrrr}
\hline\noalign{\smallskip}
\noalign{\smallskip}
Region	&	Size	&	$Structured$	& $O$ &	$DistMesh_a$	&	$O_1$ & $DistMesh_b$ & $O_2$	\\
\noalign{\smallskip}\hline\noalign{\smallskip}
dom	&	21	&	9.70E-04	&		&	4.22E-04	&		&	5.49E-04	&		\\
	&	41	&	2.26E-04	&	2.18	&	7.55E-05	&	2.57	&	1.47E-04	&	1.96	\\
	&	81	&	4.74E-05	&	2.29	&	1.04E-05	&	2.91	&	3.15E-05	&	2.27	\\
eng	&	21	&	9.69E-04	&		&	6.33E-04	&		&	6.94E-04	&		\\
	&	41	&	1.81E-04	&	2.50	&	2.47E-05	&	4.85	&	2.60E-04	&	1.47	\\
	&	81	&	4.98E-05	&	1.90	&	4.43E-05	&	-0.85	&	1.37E-04	&	0.94	\\
hab	&	21	&	1.05E-03	&		&	1.96E-04	&		&	3.55E-04	&		\\
	&	41	&	3.84E-04	&	1.51	&	2.76E-05	&	2.93	&	6.57E-05	&	2.52	\\
	&	81	&	9.54E-05	&	2.05	&	1.37E-05	&	1.03	&	2.00E-05	&	1.75	\\
m19	&	21	&	1.29E-03	&		&	2.91E-04	&		&	2.26E-03	&		\\
	&	41	&	2.72E-04	&	2.33	&	2.46E-05	&	3.70	&	1.24E-04	&	4.34	\\
	&	81	&	8.38E-05	&	1.73	&	1.36E-04	&	-2.52	&	2.87E-05	&	2.15	\\
mex	&	21	&	2.01E-03	&		&	5.55E-04	&		&	3.18E-04	&		\\
	&	41	&	3.04E-04	&	2.82	&	8.10E-05	&	2.88	&	9.42E-05	&	1.82	\\
	&	81	&	7.32E-05	&	2.09	&	1.45E-05	&	2.53	&	4.65E-06	&	4.42	\\
mic	&	21	&	8.13E-04	&		&	2.81E-05	&		&	1.10E-04	&		\\
	&	41	&	2.03E-04	&	2.07	&	7.72E-05	&	-1.51	&	8.61E-05	&	0.37	\\
	&	81	&	6.77E-05	&	1.61	&	4.39E-05	&	0.83	&	4.96E-06	&	4.19	\\
plo	&	21	&	2.36E-04	&		&	3.05E-04	&		&	1.14E-04	&		\\
	&	41	&	9.84E-05	&	1.31	&	1.23E-04	&	1.36	&	2.28E-05	&	2.41	\\
	&	81	&	2.31E-05	&	2.13	&	2.30E-05	&	2.46	&	5.04E-06	&	2.22	\\
swa	&	21	&	4.44E-04	&		&	2.44E-04	&		&	3.93E-04	&		\\
	&	41	&	1.07E-04	&	2.12	&	7.82E-05	&	1.70	&	9.46E-05	&	2.13	\\
	&	81	&	3.32E-05	&	1.72	&	1.21E-05	&	2.74	&	2.18E-05	&	2.15	\\
uch	&	21	&	1.32E-03	&		&	1.01E-03	&		&	4.97E-04	&		\\
	&	41	&	2.57E-04	&	2.45	&	2.29E-05	&	5.66	&	1.35E-04	&	1.95	\\
	&	81	&	6.26E-05	&	2.07	&	1.35E-04	&	-2.60	&	3.09E-05	&	2.17	\\
\noalign{\smallskip}\hline
\end{tabular}
\end{table}

\begin{figure}[htb]
\begin{center}
\epsfig{file=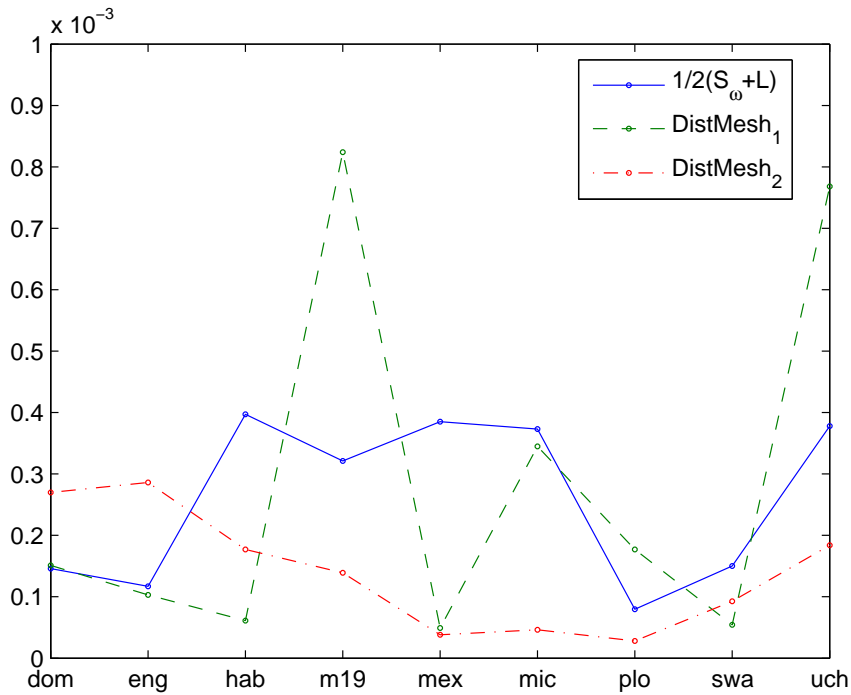, width=1.1\textwidth}
\end{center}
\caption{quadratic error for problem 1}
\label{err1}
\end{figure}

\begin{figure}[htb]
\begin{center}
\epsfig{file=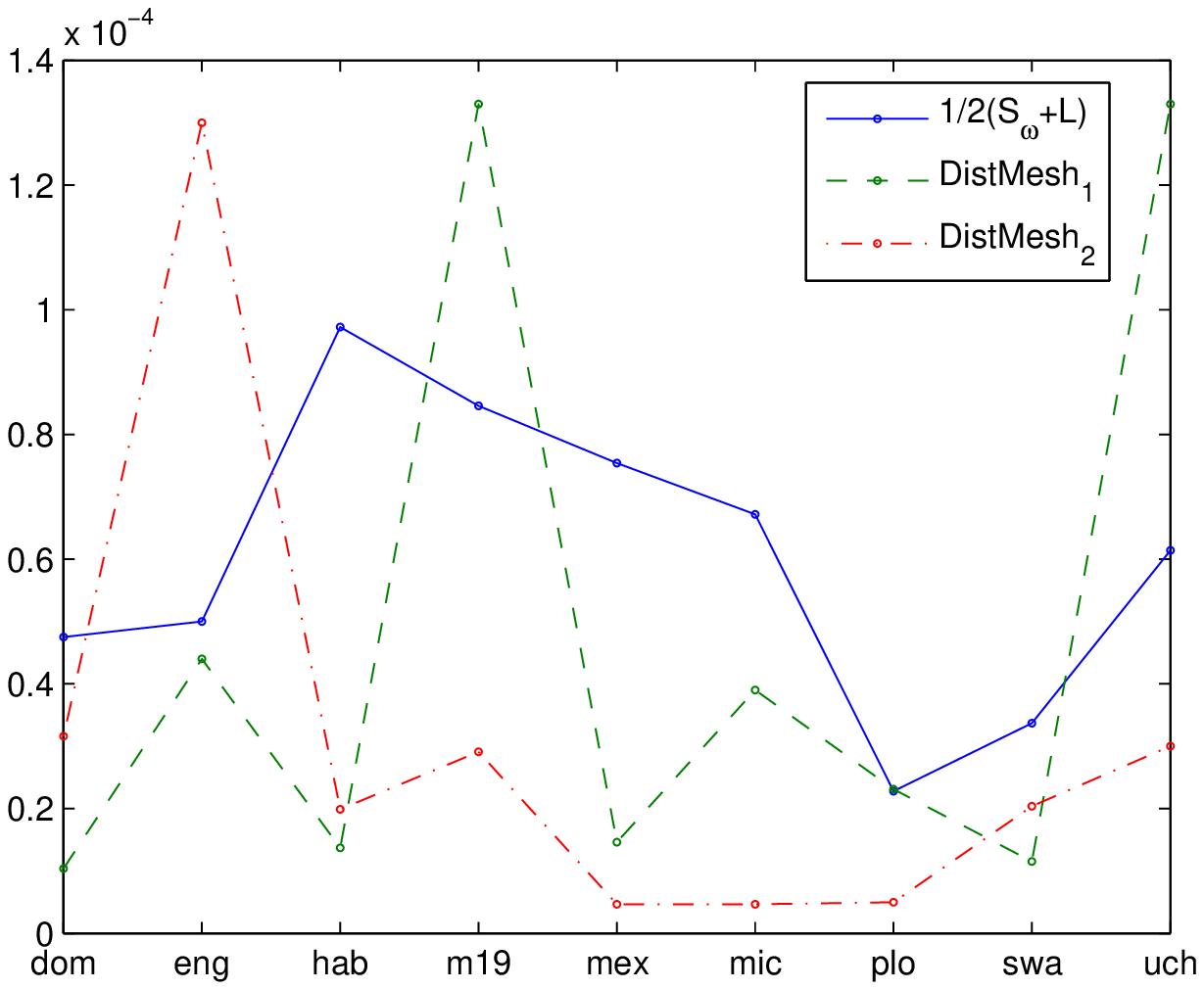, width=1.1\textwidth}
\end{center}
\caption{quadratic error for problem 2}
\label{err2}
\end{figure}

\begin{figure}[htb]
\begin{center}
\epsfig{file=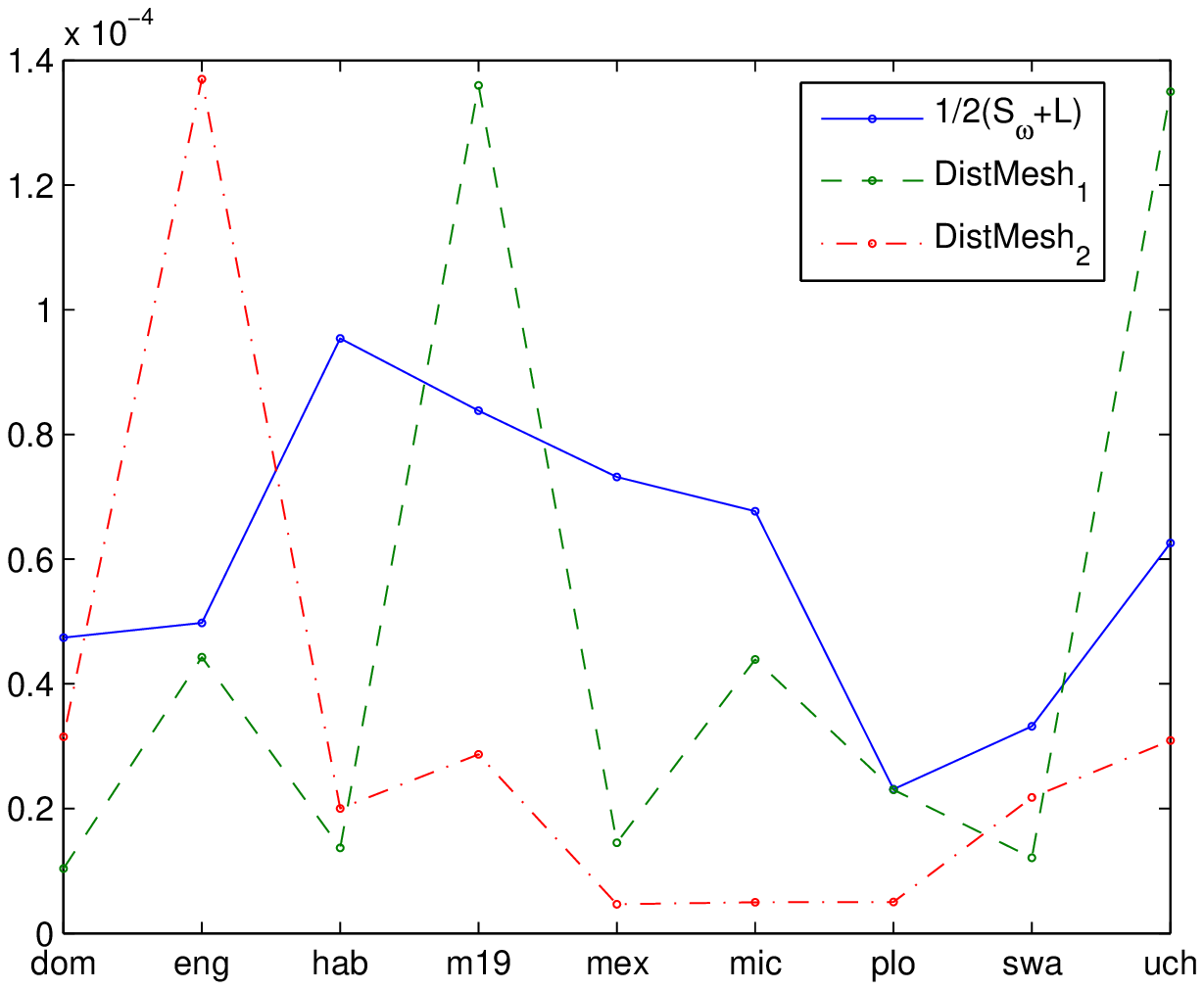, width=1.1\textwidth}
\end{center}
\caption{quadratic error for problem 3}
\label{err3}
\end{figure}

\noindent
The quadratic error for the grids with 81 points per side for the three problems is sketched in figures \ref{err1}, \ref{err2} and \ref{err3}.\\
Some conclusions can be drawn from the numerical results. The strong non convexities on some boundaries are clearly reflected in the error magnitudes and in a slight decrease of the empirical order; this can be explained in terms of the presence of elongated elements. For instance, in the $DistMesh_a$ grids, the boundary triangles must preserve the boundary shape, so the boundary nodes are kept fixed (see fig. \ref{rrr42}) and in consequence some elongated elements are produced. Nevertheless, one must also note that the very elongated cells in the structured grids have less negative effect in the solution than expected.\\
It can also be seen that in a number of problems the solution calculated with finite elements is slightly more accurate than than solution obtained with finite differences. But the former often required a larger number of unknowns, which increases the computational effort required for the calculations.\\
However, at the end, one conclusion arises: in the experiments, no method seems to be clearly superior.  The numerical solutions obtained by second order finite differences in the structured grids generated by variational methods are essentially as accurate as that obtained by second order finite elements using Delaunay-like triangulations. Having in mind this fact, an additional important issue must be discussed:  the simplicity of use of the data structures required for finite differences. Since structured grids are logically rectangular, equation (\ref{defdf}) is algorithmically as simple as a nine point discrete laplacian. This is an advantage over the more complicated data structures required for triangulations, since the numerical results do not show a significative improvement.\\
\section{Conclusions and future work}
In this paper, we selected irregular planar regions instead of rectangular ones, which are often a poor approximation to real-world domains. In addition, the side selection was quite arbitrary, and no boundary was processed in any way in order to avoid simplifying the problems. Thus, strong non convex boundaries are clearly reflected in the quadratic errors ({\em e.g. Great Britain, M19 and Ucha}), but these choices were done so because we wanted to deal with a ``hard" problem, although it must be acknowledged that boundary preprocessing is an excellent strategy to improve numerical results in general grids.\\
However, the use of irregular fixed boundaries is precisely what leads to the main conclusion: as follows from the numerical experiments, no method seems to be notably more accurate than the others; difference schemes applied on structured meshes can indeed be used to produce reliable approximations to the solution of the test problems. In other words, triangulations are not the only reliable choice for such regions, it is also possible to generate accurate numerical solutions using grids generated by variational methods in very irregular regions, and this approach has not been thoroughly studied yet. As we mentioned, this is an important issue, since differences are based on logically rectangular grids, and, in consequence, their algorithmic implementations are rather simple.\\
Our current research deals with time-dependent partial differential equations on irregular domains, and the corresponding results will be reported in a future paper.\\
\section{Acknowledgements}
We want to thank  CIC-UMSNH  ``Complejidad num\'erica y computacional
de la soluci\'on num\'erica de ecuaciones diferenciales parciales y
algunas aplicaciones-IV".
Many thanks due to the reviewers of this paper for their always valuable suggestions.

\end{document}